\documentclass{article}
\usepackage{amssymb}
\usepackage{epsfig}
\begin{document}
\title{Spectral stability of the Coulomb-Dirac Hamiltonian with anomalous magnetic moment}
%\titlerunning{Dirac Hamiltonian with anomalous magnetic moment}
\author{Hubert Kalf\thanks{Mathematisches Institut der Universit\"at, Theresienstr. 39, D-80333 M\"unchen, Germany}\ and Karl Michael Schmidt\thanks{School of Mathematics, Cardiff University, Senghennydd Road, Cardiff CF24 4YH, UK}
}
\date{}
\maketitle
\def\thet{\vartheta}
\def\ro{\varrho}
\def\eps{\varepsilon}
\def\ph{\varphi}
\def\curl{\mathop{\rm curl}}
\def\sgn{\mathop{\rm sgn}}
\def\mod{\mathop{\rm mod}}
\def\Any{\mathord{\cdot}}
\def\R{\mathbb{R}}
\def\Z{\mathbb{Z}}
\def\N{\mathbb{N}}
\def\modul#1{\mathopen{\vert}#1\mathclose{\vert}}
{\it \qquad Dedicated to A. Schneider on the occasion of his 70th birthday\/}
\medskip
\begin{abstract}
We show that the point spectrum of the standard Coulomb-Dirac operator $H_0$ is the limit of
the point spectrum of the Dirac operator with anomalous magnetic moment $H_a$ as the anomaly
parameter tends to $0$.
For negative angular momentum quantum number $\kappa$, this holds for all Coulomb coupling
constants $c$ for which $H_0$ has a distinguished self-adjoint realisation.
For positive $\kappa$, however, there are some exceptional values for $c$, and in general an
index shift between the eigenvalues of $H_0$ and the limits of eigenvalues of $H_a$ appears,
accompanied with additional oscillations of the eigenfunctions of $H_a$ very close to the
origin.
\end{abstract}
\section{Introduction}
\label{sec:intro}

In a 1930 letter to O. Klein (\cite{P2} letter 261) and in a survey to be given at the 8th Solvay
Congress, planned for October 1939, Pauli suggested to describe the motion of a particle
with rest mass $m > 0$, charge $e$, spin $\hbar/2$ and magnetic moment $(1 + a) \mu_B$
($\mu_B = {e \hbar \over 2 m c}$ the Bohr magneton) in an electric field $-\nabla \Phi$ and
a magnetic field $B = \curl A$ by means of the operator
$$
  H = c\, \alpha\cdot(p - {e \over c} A) + m c^2 \beta + e \Phi - a \mu_B(i \alpha\cdot\nabla\Phi
      + \sigma\cdot B).
$$
A revised version of this review appeared in \cite{P1} (the relevant equation is (91)); the
original manuscript was not published until 1993 in \cite{P3} pp. 827--901.

If $\Phi(x) = V(\modul{x})$ and $A = 0$, self-adjoint realisations of $H$ in $L^2(\R^3)^4$ are
unitarily equivalent to the orthogonal sum of self-adjoint realisations of
$$
  H_a = \sigma_2 p + {m c \over \hbar} \sigma_3 + \left({\kappa \over r}
      + {a \mu_B \over \hbar c} V'(r)\right)
        \sigma_1 + {e \over \hbar c} V(r)
$$
in $L^2((0, \infty))^2$, where $\kappa \in \Z\setminus\{0\}$ is the angular momentum quantum number.

The first mathematical treatment of the operators $H$ and $H_a$ is due to Behncke \cite{B1,B2,B3}.
He showed that $H_a$ has a unique self-adjoint realisation if $a \neq 0$
for a very large class of potentials $V$, including the Coulomb potential $V(r) = -{Z e \over r}$
for all values of the coupling constant $Z e$
(for alternative proofs see \cite{GST} and \cite{AKS}).
This is in marked contrast to the case $a = 0$ where it is well known that $H_0$ is essentially self-adjoint
on its minimal domain if and only if $({Z e^2 \over \hbar c})^2 \le \kappa^2 - {1 \over 4}$.
For larger values
of $Z$ the singular end-point $0$ is in the limit-circle
case; but as long as $({Z e^2 \over \hbar c})^2 < \kappa^2$ one still has a distinguished
self-adjoint realisation
of $H_0$ defined by the requirement that functions in the domain behave like the principal
solution of the eigenvalue equation of $H_0$ at $0$.

The location of the essential spectrum of $H_a$ is
comparatively easy to determine and for a wide range of
potentials, notably the Coulomb potential, coincides
with that of $H_0$ \cite{B1,B2,B3}.

In the following, we normalise constants and write the
Coulomb Hamiltonian as
$$
  H_a = -i \sigma_2 {d \over dr} + \sigma_3 + \left({\kappa \over r} + {a \over r^2}\right)
        \sigma_1 + {c \over r},
$$
assuming $a, c < 0$.
(The cases of positive $c$ and/or $a$ can be reduced to this situation by means of suitable
unitary transformations.)
The discrete spectrum of $H_0$ accumulates at the right end-point
of the gap $(-1, 1)$ in the essential spectrum. Since
$C_0^\infty((0, \infty))^2$ is a common core for $H_0$ and $H_a$ if (and only if)
$c^2 \le \kappa^2 - {1 \over 4}$, $H_a$ converges to $H_0$ in the strong resolvent
sense as $a \rightarrow 0$ (\cite{RS} Thm VIII.25(a)), and
as a consequence, the spectrum of $H_a$ cannot expand in the limit $a \rightarrow 0$ for this range
of the parameters $c, \kappa$. However, it is reasonable to expect (and has been used as
a basis for a perturbative calculation of the eigenvalues of $H_a$) that the point
spectrum is stable in the limit $a \rightarrow 0$ in the sense that the eigenvalues
of $H_a$ converge to those of $H_0$, and each eigenvalue of $H_0$ is the limit of
exactly one eigenvalue branch of $H_a$.
Decoupling the eigenvalue equation of $H_a$ and using a comparison theorem for principal
and non-principal solutions of second-order equations, Behncke \cite{B3} proved this
stability for
$\kappa^2 > c^2 + ({3\over 2})^2$ if $\kappa < 0$,
and
$\kappa > c^2 + {5 \over 2}$ if $\kappa > 0$.
He conjectured that $\kappa^2 > c^2 + {5 \over 2}$ might be sufficient in the latter case.
(Farther-reaching conjectures are to be found in \cite{Th} p. 218 seq.)

\smallskip
In the present paper, we study the convergence of the point spectrum of $H_a$ as $a$
tends to $0$, for the whole parameter range for which a distinguished realisation of $H_0$
exists, i.e. for $\kappa^2 - c^2 > 0$.
We find a surprising qualitative difference in the limiting behaviour depending on the
sign of $\kappa$. Indeed, for negative $\kappa$, the eigenvalues of $H_0$ are exactly
the limits of eigenvalues of $H_a$ for all values of $c$.
For positive $\kappa$, however, there are (finitely or infinitely many) exceptional values
$c_0 > c_1 > \dots$ in $(-\kappa, 0)$; for $c \in (c_m, c_{m-1})$, the eigenvalues of
$H_0$ are still the limits of the eigenvalues of $H_a$, but with a shift of size $m$ in
the eigenvalue numbers. This shift is reflected in the appearance of $m$ additional
oscillations of the corresponding eigenfunction of $H_a$, compared to that of $H_0$, very
close to the origin.
It seems a delicate question to decide whether the number of exceptional values $c_m$ is
finite or infinite; in any case it grows beyond all bounds with increasing $\kappa$.

More precisely, we have the following results.

\medskip
\noindent
{\bf Theorem 1.1 %\begin{theorem}
(Spectral convergence and stability for negative~$\kappa$)}
%\label{th:11}

{\it
\noindent
Let $\kappa < 0$, $c \in (\kappa, 0)$, and let
$\lambda_0$ [not] be an eigenvalue of the Coulomb-Dirac Hamiltonian 
$$
  H_0 = -i\sigma_2 {d \over dr} + \sigma_3 + {\kappa \over r} \sigma_1 + {c \over r}.
$$
Let $0 < \varepsilon < \mathop{\rm dist}(\lambda_0, \sigma(H_0) \setminus
 \{\lambda_0\})/2$.
Then for $a < 0$ with sufficiently small $\modul a$ the
Hamiltonian with anomalous magnetic moment
$$
  H_a = -i\sigma_2 {d \over dr} + \sigma_3 + \left({\kappa \over r} + {a \over r^2}\right)
     \sigma_1 +   {c \over r}
$$
has exactly one [no] eigenvalue $\lambda_a$ in
$(\lambda_0 - \varepsilon, \lambda_0 + \varepsilon)$.

}
%\end{theorem}

\medskip
\noindent
{\bf Theorem 1.2 %\begin{theorem}
(Spectral convergence and stability for positive~$\kappa$)}
%\label{th:12}

{\it
\noindent
Let $\kappa > 0$; then there are at least $[(\kappa \log 4)/\pi - 1]$ values
$0 > c_0 > c_1 > \dots > -\kappa$, which can only accumulate at $-\kappa$, such that the following
holds.

Let $c \in (-\kappa, 0) \setminus \{c_0, c_1, \dots\}$, and let
$\lambda_0$ [not] be an eigenvalue of $H_0$.
Let $0 < \varepsilon < \mathop{\rm dist}(\lambda_0, \sigma(H_0) \setminus
 \{\lambda_0\})/2$.
Then for $a < 0$ with sufficiently small $\modul a$,
$H_a$ has exactly one [no] eigenvalue $\lambda_a$ in
$(\lambda_0 - \varepsilon, \lambda_0 + \varepsilon)$.

}
%\end{theorem}

\medskip
\noindent
Here $[x] := \sup\{m \in \Z \mathrel{|} m \le x\}$ $(x \in \R)$ denotes the
Gau\ss\ bracket.

The proof of these theorems is based on oscillation theory, in particular on an asymptotic
study of the behaviour of the Pr\"ufer angle of solutions of the eigenvalue equation for
$H_a$ (i.e. the solutions of equation (\ref{eq:31}) below) as $a$ tends to $0$.
After rescaling $\varrho = r/\modul{a}$, the mass term and spectral parameter are lower
order terms in the limit and can be omitted in order to obtain an overview of solutions
in the asymptotic regime near the origin.

The direction field of the resulting simplified equation (\ref{eq:21}) (where $k = \modul\kappa$
and $\alpha = -\sgn\kappa$), and hence the qualitative behaviour of its solutions, shows a
fundamental difference depending on the sign of $\alpha$.

\begin{figure}[t] 
  \begin{center} 
    \epsfig{file=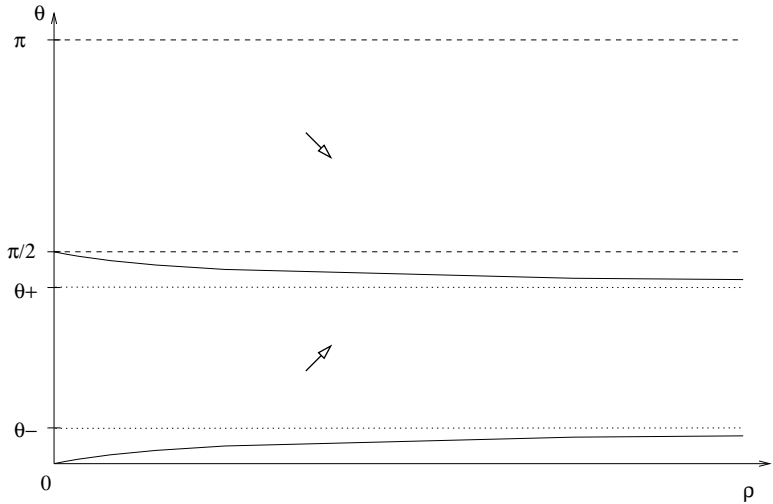,width=4in} 
    \caption{\label{fig:kapneg}Zones in the directional field of the simplified equation for $\alpha = 1$.
The arrows indicate the sign of the right-hand side of (\ref{eq:21}), the solid lines represent its zeros.} 
  \end{center} 
\end{figure}

\medskip
For negative $\kappa$ ($\alpha = 1$), the $(\ro, \thet)$ plane is divided into essentially
horizontal strips in which the right-hand side of (\ref{eq:21}) is alternatingly positive and
negative (cf. Fig. \ref{fig:kapneg}); as a consequence, the distinguished angle $\thet_0$ which corresponds
to an $L^2(0, \Any)$ solution of the eigenvalue equation for $H_a$, cannot change
by more than $\pi$, and eventually tends towards an asymptotic value $\thet_+(c)$
(Prop. 2.3 a) which turns out to be the limiting angle at $0$ of the principal solution
of the eigenvalue equation for $H_0$ as well (Lemma 3.1 b).
On the original $r$ scale, the convergence becomes faster as $a \rightarrow 0$.

A stability argument (Lemma 3.2, Prop. 3.3) then shows that for a certain point $R > 0$
the influence of the previously neglected mass and spectral parameter terms can be controlled
on $(0, R)$, and that the solution of the full Pr\"ufer equation (\ref{eq:31}) converges at that
point to the Pr\"ufer angle of the principal solution of the equation for $a = 0$.

Theorem 1.1 follows in view of the uniformity of this convergence with respect to the
spectral parameter, and the fact that the presence of the anomalous magnetic moment does not
significantly affect the behaviour of the solutions at $\infty$.

\begin{figure}[t] 
  \begin{center} 
    \epsfig{file=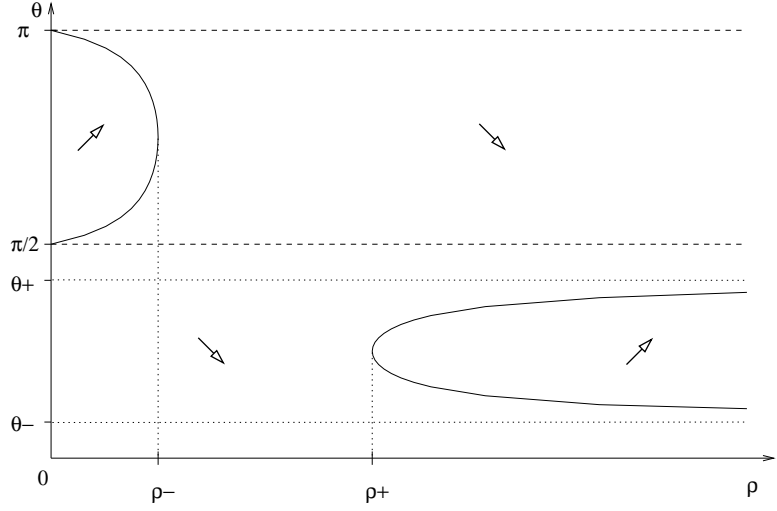,width=4in} 
    \caption{\label{fig:kappos}Zones in the directional field of the simplified equation for $\alpha = -1$.
The arrows indicate the sign of the right-hand side of (\ref{eq:21}), the solid lines represent its zeros.} 
  \end{center} 
\end{figure}

\medskip
A curious phenomenon occurs, however, in the case of positive $\kappa$ ($\alpha = -1$).
For $\varrho$ close to $0$ and for large $\varrho$, one again has $\thet$-regions of
opposite sign of the right-hand side of (\ref{eq:21}), and hence of essential confinement of
the solutions, and for large $\varrho$, the distinguished angle $\thet_0$ generically
tends to the limit $\thet_+(c) \mod \pi$.
In contrast to the previous situation, there is now a $\varrho$-interval $(\ro_-(c),
\ro_+(c))$ on which the right-hand side of (\ref{eq:21}) has no zeros and is strictly negative
(cf. Fig. \ref{fig:kappos}).

Depending on $c$, the size of this gap increases, vanishing as $c \rightarrow 0$ and
becoming infinite as $c \rightarrow -\kappa$.
Moreover, one can show that the angle $\thet_0$ changes by several multiples of $\pi$
in this interval if $-c$ is large enough (Prop. 2.3 b), and hence will eventually
converge to $\thet_+(c) - m\pi$ for some $m \in \N_0$.
Thus the corresponding $L^2(0, \Any)$ solution of the eigenvalue equation for $H_a$
will perform $m$ oscillations which are absent in the principal solution of the
eigenvalue equation for $H_0$, to which it however converges in phase at the point
$R > 0$.

Theorem 1.2 then follows by way of the same stability argument as before (Lemma 3.2,
Prop. 3.3).

As $c$ crosses the exceptional value, a transition of the asymptotic limit from
$\thet_+(c) - m\pi$ to $\thet_+(c) - (m+1)\pi$ takes place.
At $c = c_m$, the distinguished solution $\thet_0$ of the simplified equation (\ref{eq:21})
approaches an unstable limit as $\ro\rightarrow\infty$.
The asymptotics of the solutions of the full Pr\"ufer equation (\ref{eq:31}) appear to be a
rather delicate matter in these cases, and the limiting behaviour of $H_a$ with
$c = c_m$ remains an interesting open question.

The basic analytical tool of the present paper is the study of the direction field of
Pr\"ufer and Riccati type ordinary differential equations near a singularity
(for related earlier but simpler results see \cite{K} and the references therein).
The underlying comparison techniques are outlined in the Appendix.

\section{The simplified equation}
\label{sec:simpeq}

In this section we consider the scaled and simplified equation
\begin{equation}
\label{eq:21}
  \varrho\,\vartheta' = c + (k + \alpha/\varrho)\,\sin 2\vartheta
\end{equation}
(with $k > 0$, $c \in (-k, 0)$ and $\alpha \in \{-1, 1\}$)
which arises from the Pr\"ufer equation (\ref{eq:31}) equivalent to the eigenvalue
equation for $H_a$ by omitting the $O(1)$ $(\varrho\rightarrow\infty)$ terms
after rescaling $r = \modul{a}\varrho$, which eliminates $\modul{a}$.

The asymptotic zeros $\thet_\pm(c)$ of the right-hand side of equation (\ref{eq:21}) satisfy
$$
\begin{array}{c}
  0 < \vartheta_-(c) < \pi/4 < \vartheta_+(c) < \pi/2,
\cr\cr
  \sin 2 \vartheta_\pm(c) = -c/k,
\quad
  \tan \vartheta_\pm(c) = {k \pm \sqrt{k^2 - c^2} \over -c}.
\cr
\end{array}
$$
Moreover, in case $\alpha = -1$ there are the exceptional points
$\varrho_\pm(c) = 1/(k \pm c)$,
between which the right-hand side of (\ref{eq:21}) is strictly negative;
for later convenience we define $\varrho_+(c) = 0$ if $\alpha = 1$.

We first show that (\ref{eq:21}) has a stable and an unstable asymptotic critical point
at both singular end-points $0$ and $\infty$ (Lemma 2.1, 2.2).
The distinguished (unstable) solution at $0$, $\thet_0$, converges to the stable
limit at $\infty$ for all $c$ if $\alpha = 1$; if $\alpha = -1$, one has the
convergence with a shift by an integer multiple of $\pi$ unless $c$ is one of a
sequence of exceptional values (Prop. 2.3).

\medskip
\noindent
{\bf Lemma 2.1 %\begin{lemma}
(The existence of a distinguished solution at $\infty$)}
\label{lem:21}

{\it
\noindent
For each $c \in (-k, 0)$, (\ref{eq:21}) has a unique solution
$\thet_\infty(\Any, c)$ such that
$$
  \lim_{\varrho\rightarrow\infty} \thet_\infty(\varrho, c)
   = \vartheta_-(c).
$$
All other solutions $\vartheta(\Any, c)$ satisfy either
$
  \vartheta(\varrho, c) = \thet_\infty(\varrho, c) + m \pi
\qquad
  (\varrho > 0)
$
for some $m \in \Z$,
or else
$
  \lim_{\varrho\rightarrow\infty} \vartheta(\varrho, c)
   = \vartheta_+(c) \mod \pi.
$
Furthermore, for each $\ro > 0$, $\thet_\infty(\ro, \Any)$ is continuous and
strictly decreasing.
\/} %\end{lemma}

\medskip
\noindent
The {\it proof\/} of this lemma, based on an asymptotic study of the direction
field of (\ref{eq:21}), can be found in the Appendix.
Similarly, one can prove

\medskip
\noindent
{\bf Lemma 2.2 %\begin{lemma}
(The existence of a distinguished solution at $0$)}
\label{lem:22}

{\it
\noindent
For each $c \in (-k, 0)$ there is a unique solution $\thet_0(\Any, c)$
such that
$$
  \lim_{\ro\to 0} \thet_0(\ro, c) = \cases{\pi & if $\alpha = -1$, \cr
                                           \pi/2 & if $\alpha = 1$. \cr}
$$
All other solutions are either shifts of $\thet_0$ by an integer
multiple of $\pi$, or have
$$
  \lim_{\ro\to 0} \thet(\ro, c)
   = \cases{\pi/2 \mod \pi& if $\alpha = -1$, \cr
            0 \mod \pi & if $\alpha = 1$. \cr}
$$
For fixed $\hat\ro > 0$, $\thet_0(\hat\ro, \Any)$ is continuous
non-decreasing.
\/} %\end{lemma}

\medskip
\noindent
{\bf Proposition 2.3} %\begin{proposition}
\label{prop:23}

{\it
\noindent
a)
If $\alpha = 1$, we have
$\lim\limits_{\ro\to\infty} \thet_0(\ro, c) = \thet_+(c)$
$(c \in (-k, 0))$.

\noindent
b)
If $\alpha = -1$, there are at least
$[(k \log 4)/\pi - 1]$ values
$0 > c_0 > c_1 > c_2 > \dots > -k$
(accumulating at $-k$ if infinitely many)
such that
$$
  \lim_{\ro\to\infty} \thet_0(\ro, c) = \thet_+(c) - m \pi
\qquad
  (c \in {\cal C}_m),
$$
and
$
  \lim\limits_{\ro\to\infty} \thet_0(\ro, c_m) = \thet_-(c_m) - m \pi.
$

Here ${\cal C}_m := (c_m, c_{m-1})$ (where $c_{-1} := 0$) for all
$m\in \N_0$ for which $c_m$ exists; if there is a minimal $c_{m_{\rm max}} > -k$,
define ${\cal C}_{m_{\rm max} + 1} := (-k, c_{m_{\rm max}})$.
\/} %\end{proposition}

\medskip
\noindent
\textit{Remark.}
For later convenience, we also define ${\cal C}_0 := (-k, 0)$ in the case
$\alpha = 1$.

\medskip
\noindent
{\it Proof.\/} %\begin{proof}
a)
The interval $[\pi/4, \pi]$ is stable for (\ref{eq:21}) by Lemma 4.2 (cf. the
Appendix), so $\thet_0(\Any, c)$ (which by Lemma 2.2 is close to $\pi/2$
for small $\ro$) cannot tend to $\thet_-(c) \mod \pi$.
The assertion follows by Lemma 2.1.

\noindent
b)
In the limit $c \rightarrow 0$, we have, using
$\thet'(\ro) \ge (1/\ro) (c - \modul{k - 1/\ro})$,
$$
\begin{array}{rl}
  &\thet_0(\ro_+(c), c) - \thet_0(\ro_-(c), c)
  \ge \int_{\ro_-(c)}^{1/k} {1 \over \ro}\,(c - {1 \over \ro} + k)\,d\ro
   +  \int_{1/k}^{\ro_+(c)} {1 \over \ro}\,(c + {1 \over \ro} - k)\,d\ro
\cr\cr
  &\qquad\qquad = k - {1 \over \ro_-(c)} - (c + k) \log(k\ro_-(c))
   - {1 \over \ro_+(c)} + k + (c - k) \log(k\ro_+(c))
\cr\cr
  &\qquad\qquad = (c - k) \log {k \over k + c} - (c + k) \log {k \over k - c}
  \to 0.
\cr
\end{array}
$$
As a consequence, $\lim\limits_{\ro\to 0} \thet_0(\ro, c) = \thet_+(c)$
for $c < 0$ sufficiently close to $0$.

On the other hand, noting that
$\thet'(\ro) \le (1/\ro) (c + \modul{k - 1/\ro})$ and
$\thet_\infty(\ro_+(c), c) \in [0, \pi/4]$, we have
$$
\begin{array}{rl}
  &\thet_0(1/k, c) - \thet_\infty(1/k, c)
  \le \thet_0(\ro_-(c), c)
   + \int_{\ro_-(c)}^{1/k} {1 \over \ro}\,(c + {1 \over \ro} - k)\,d\ro
\cr\cr
  &\qquad\qquad\quad
   - \thet_\infty(\ro_+(c), c)
   + \int_{1/k}^{\ro_+(c)} {1 \over \ro}\,(c - {1 \over \ro} + k)\,d\ro
\cr\cr
  &\qquad\le \pi - k + {1 \over \ro_-(c)} + (k - c) \log(k\ro_-(c))
   + {1 \over \ro_+(c)} - k + (k + c) \log(k\ro_+(c))
\cr\cr
  &\qquad\to \pi - k \log 4
\qquad
  (c \to -k).
\cr
\end{array}
$$
The assertion follows, as $\thet_0(1/k, \Any)$ and $\thet_\infty(1/k, \Any)$
are continuous and monotone non-decreasing and decreasing, resp.
\hfill $\square$ %\qed\end{proof}

\medskip
\noindent
\textit{Remark.}
The first part of the proof of Proposition 2.3 b) yields a quantitative
estimate for the first exceptional value $c_0$.
Indeed, using Lemma 4.1 one can see that
$\thet_0(\ro_-(c_0), c_0) \ge {3 \pi \over 4}$
and
$\thet_0(\ro_+(c_0), c_0) \le {\pi \over 4}$, so
$$
  {\pi \over 2} \le (c_0 - k) \log{k + c_0 \over k} - (c_0 + k) \log{k - c_0 \over k}.
$$
Setting $x := -c_0 / k \in (0, 1)$ and observing that
$$
\begin{array}{rl}
  x \log {1 + x \over 1 - x} - \log (1 - x^2)
  &= \sum_{j=1}^\infty x^{2j} \left({1 \over j} + {1 \over {j - {1 \over 2}}}\right)
\cr\cr
  &= x^2 \left(3 + x^2 \sum_{j=0}^\infty x^{2j}
    \left({1 \over j + 2} + {1 \over {j + {3 \over 2}}}\right)\right)
\cr\cr
  &\le x^2 \left(3 + {x^2 \over 1 - x^2}\,{7 \over 6}\right)
  = {x^2 (18 - 11 x^2) \over 6 (1 - x^2)},
\cr
\end{array}
$$
we find
$11 x^4 - (18 + 3\pi/k) x^2 + 3\pi/k \le 0$, and thus
$$
  c_0^2 \ge {1 \over 22}\,(18 k^2 + 3 \pi k - \sqrt{(18 k^2 + 3 \pi k)^2 - 132 \pi k^3}).
$$
For $k = 1$, this gives the bound $c_0 \le -0.64157$.

Hence, one has convergence and stability of the eigenvalues by Theorem 1.2 for all
$\kappa \in \Z \setminus \{0\}$ at least for nuclear charge number $Z \le 87$.

\section{The convergence of the original equation for $a \to 0$, $c \in {\cal C}_m$}
\label{sec:conv}

Throughout this section, we fix $\alpha, \mu \in \{-1, 1\}$
and $k > 0$, and assume that $c \in {\cal C}_m$ for some admissible $m \in \N_0$.

Consider the Pr\"ufer equation
\begin{equation}
\label{eq:31}
  \Theta' = {c \over r} + ({k \over r} + {\alpha\modul{a} \over r^2}) \sin 2 \Theta
  + \mu \cos 2 \Theta - \lambda
\end{equation}
and the corresponding equation for $a = 0$,
\begin{equation}
\label{eq:32}
  X' = {c \over r} + {k \over r} \sin 2 X
  + \mu \cos 2 X - \lambda.
\end{equation}

We shall study the following distinguished solutions, whose existence and properties
can be proved in analogy to Lemma 2.1.

\medskip
\noindent
{\bf Lemma 3.1} %\begin{lemma}
\label{lem:31}

{\it
\noindent
a)
There is exactly one solution $\Theta_0$ of (\ref{eq:31}) with
$$
  \lim_{r\rightarrow 0} \Theta_0(r, c, a, \lambda)
   = \cases{ \pi &if $\alpha = -1$ \cr
             \pi/2 &if $\alpha = 1$, \cr}
$$
all other solutions either being shifts of $\Theta_0$ by an integer multiple
of $\pi$, or having
$$
  \lim_{r\rightarrow 0} \Theta(r, c, a, \lambda)
   = \cases{\pi/2\mod\pi &if $\alpha = -1$\cr                                      
            0\mod\pi &if $\alpha = 1$. \cr}
$$
For fixed $r, c$ and $a$, $\Theta_0(r, c, a, \Any)$ is continuous decreasing.

\noindent
b)
There is exactly one solution $X_0$ of (\ref{eq:32}) with
$$
  \lim_{r\rightarrow 0} X_0(r, c, \lambda) = \thet_+(c),
$$
all other solutions either being shifts of $X_0$ by an integer multiple
of $\pi$, or having
$$
  \lim_{r\rightarrow 0} X(r, c, \lambda) = \thet_-(c) \mod \pi.
$$
For fixed $r$ and $c$, $X_0(r, c, \Any)$ is continuous
decreasing.
\/} %\end{lemma}

\medskip
\noindent
Now we show that the solutions $X_0$ and $\Theta_0$ become asymptotically close $\mod \pi$
at some point for small $a$; this is a consquence of the convergence of the solution $\thet_0$
(which is close to $\Theta_0$) of the simplified equation (Prop. 2.3).

\medskip
\noindent
{\bf Lemma 3.2} %\begin{lemma}
\label{lem:32}
{\it
For each $\eps > 0$ there are $r_0(\eps) \in (0, \eps]$ and
$a_0 > 0$ such that for all $\lambda \in [-1, 1]$ and
$-a_0 < a < 0$
$$
\begin{array}{rl}
  \modul{\Theta_0(r_0(\eps), c, a, \lambda) + m\pi - \thet_+(c)} &< {2\eps\over
3}, \cr\cr
  \modul{X_0(r_0(\eps), c, \lambda) - \thet_+(c)} &< {\eps\over 3},
\cr
\end{array}
$$
and consequently
$
  \modul{\Theta_0(r_0(\eps), c, a, \lambda) + m\pi
    - X_0(r_0(\eps), c, \lambda)} < \eps.
$
\/} %\end{lemma}

\medskip
\noindent
{\it Proof.\/} %\begin{proof}
Let $\gamma > 0$ be so small that
$[c - \gamma, c + \gamma] \subset {\cal C}_m$.
Then there is $\delta \in (0, \gamma)$ such that
$\modul{\thet_+(c \pm \delta) - \thet_+(c)} < \eps/3$.
By Proposition 2.3 there is $\ro_0 > 0$ such that
$\modul{\thet_0(\ro, c \pm \delta) - \thet_+(c \pm \delta) + m \pi}
 < \eps/3$
$(\ro \ge \ro_0)$,
so the function $r \mapsto \thet_0(r/\modul{a}, c)$
(which is a solution of the simplified equation
$$
  \tilde\Theta'(r)
   = {c \over r} + ({k \over r} + {\alpha\modul{a} \over r^2}) \sin 2 \tilde\Theta)
$$
satisfies
$
  \modul{\thet_0(r/\modul{a}, c \pm \delta) - \thet_+(c) + m \pi}
  < 2 \eps/3
\qquad
  (\modul{a} \le {r \over \ro_0}, r > 0).
$
Now let $r_0(\eps) := \min \{\delta/2, \eps\}$.
Estimating in (\ref{eq:32}), (\ref{eq:31})
$$
  \modul{\mu \cos 2 X - \lambda}, \modul{\mu \cos 2 \Theta - \lambda}
  \le 2 \le {\delta \over r}
\qquad
  (r \le r_0(\eps)),
$$
Lemma 4.1 yields the bounds
$$
  \thet_0(r/\modul{a}, c-\delta) \le \Theta_0(r, c, a)
  \le \thet_0(r/\modul{a}, c+\delta)
\qquad
  (0 < r \le r_0(\eps))
$$
and
$$
  \thet_+(c-\delta) \le X_0(r, c) \le \thet_+(c+\delta)
\qquad
  (0 < r \le r_0(\eps)),
$$
and the assertion follows.
\hfill $\square$ %\qed\end{proof}

\medskip
\noindent
The preceding lemma shows that the solutions converge to each other; however,
the point of comparison $r_0(\eps)$ depends on $\eps$ and tends
to $0$ rather rapidly.
We now show that the convergence remains stable, and hence also holds
at a certain fixed point $R$.

\medskip
\noindent
{\bf Proposition 3.3} %\begin{proposition}
\label{prop:33}
{\it
There is $R > 0$ such that
$
  \lim\limits_{a\to 0} \Theta_0(R, c, a, \lambda) = X_0(R, c, \lambda)
$
uniformly w.r.t. $\lambda \in [-1, 1]$.
\/} %\end{proposition}

\medskip
\noindent
{\it Proof.\/} %\begin{proof}
Consider the Riccati equations for $z(r) := \tan X_0(r, c)$ and
$y(r, a) := \Theta_0(r, c, a)$,
\begin{equation}
\label{eq:33}
 r\,z' = (c - (\lambda + 1) r) z^2 + 2 k z + c - (\lambda - 1) r
\end{equation}
and
\begin{equation}
\label{eq:34}
 r\,y' = (c - (\lambda + 1) r) z^2 + 2 (k+\alpha\modul{a}/r) z + c - (\lambda - 1) r.
\end{equation}
Let
${\displaystyle
  y_+(c) := \tan \thet_+(c) = {k + \sqrt{k^2 - c^2} \over -c}.
}$
Choose $0 < d < \sqrt{k^2 - c^2}/(3 \modul{c})$ and $0 < R \le d$
so small that
$$
  {3 \over d}
   \left\vert\modul{\lambda + 1} (y_+(c) + d)^2 + \modul{\lambda - 1}\right\vert R,
  6 \modul{\lambda + 1} (y_+(c) + d) R
  < \sqrt{k^2 - c^2}
$$
for all $\lambda \in [-1, 1]$.
For any $\hat r \in (0, R)$, set
$$
  a_1(\hat r) := {\sqrt{k^2 - c^2} \over 6 ({y_+(c) \over d} + 1)} \hat r.
$$

Now let $\eps > 0$, $\eps < R \le d$.
By Lemma 3.2, there is $r_0(\eps) < \eps$ and $a_0 > 0$ such that
$$
  \modul{z(r_0(\eps)) - y_+(c)} < \eps/2,
\quad
  \modul{y(r_0(\eps), a) - y_+(c)} < \eps/2
\quad
  (0 < \modul{a} < a_0).
$$
On $[r_0(\eps), R]$, the interval $[y_+(c) - d, y_+(c) + d]$ is stable
for (\ref{eq:33}) and (\ref{eq:34}) by Lemma 4.2
if $\modul{a} < a_1(r_0(\eps))$, since the
right-hand side for $y_+(c) \pm d$ takes the value
$$
\begin{array}{rl}
  &(c - (\lambda + 1) r) (y_+(c) \pm d)^2 + 2 (k + \alpha\modul{a}/r) (y_+(c) \pm d)
  + c - (\lambda - 1) r
\cr\cr
  &\ = c d^2 \pm 2 (k + y_+(c) c) d + 2 {\alpha\modul{a} \over r} (y_+(c) \pm d)
  - (\lambda + 1) r (y_+(c) \pm d)^2 - (\lambda - 1) r
\cr\cr
  &\ \mp d \left[ 2 \sqrt{k^2 - c^2} \pm c d
   \pm 2 {\alpha\modul{a} \over r} ({y_+(c) \over d} \pm 1)
   \mp ((\lambda + 1) (y_+(c) \pm d)^2 - \lambda - 1) {r \over d} \right],
\cr
\end{array}
$$
and the factor in square brackets is not less than $\sqrt{k^2 - c^2}$.

Hence
$\modul{y(r, a) - y_+(c)}, \modul{z(r) - y_+(c)} \le d$
for all $r \in [r_0(\eps), R]$ if $0 < \modul{a} < \min\{a_0, a_1(r_0(\eps))\}$.

In the differential equation for the difference
$x(r, a) := y(r, a) - z(r)$,
\begin{equation}
\label{eq:35}
  r\,x' = (c - (\lambda + 1) r) (y + z) + 2k)\,x - {2 \alpha\modul{a} \over r} y,
\end{equation}
the factor of $x$ on the right-hand side can be estimated
$$
\begin{array}{rl}
  c (y(r, a) + z(r)) &+ 2k - (\lambda + 1) r (y + z)
\cr\cr
  &\le -2(\sqrt{k^2 - c^2} + cd) + 2 R \modul{\lambda + 1} (y_+(c) + d)
\cr\cr
  &\le - \sqrt{k^2 - c^2}.
\cr
\end{array}
$$
Hence, if $a$ additionally satisfies
$$
  \modul{a} < {\sqrt{k^2 - c^2} \,r_0(\eps) \over 2 (y_+(c) + d)}\,\eps,
$$
then the interval $[-\eps, \eps]$ is stable for (\ref{eq:35}) on $[r_0(\eps), R]$
by Lemma 4.2,
and because of $\modul{x(r_0(\eps), a)} < \eps$ we find
$\modul{y(R, a) - z(R)} < \eps$.
\hfill $\square$ %\qed\end{proof}

\medskip
\noindent
For the proof of the main theorems, we need the following distinguished
solution of (\ref{eq:32}) at infinity, whose existence and properties can be proven
along the lines of Lemma 2.1. Let $X_\pm(\lambda) \in [0, \pi]$ such that
$\cos 2 X_\pm(\lambda) = \mu \lambda$ and
$\sin 2 X_\pm(\lambda) = \pm \mu \sqrt{1 - \lambda^2}$ $(\lambda \in [-1, 1])$.

\medskip
\noindent
{\bf Lemma 3.4} %\begin{lemma}
\label{lem:34}
{\it
For $\lambda \in [-1, 1]$, there is exactly one solution $X_\infty$ of (\ref{eq:32})
with $$
  \lim_{r\rightarrow \infty} X_\infty(r, c, \lambda) = X_-(\lambda),
$$
all other solutions either being shifts of $X_0$ by an integer multiple
of $\pi$, or having
$$
  \lim_{r\rightarrow \infty} X(r, c, \lambda) = X_+(\lambda) \mod \pi.
$$
For fixed $r$ and $c$, $X_\infty(r, c, \Any)$ is continuous
increasing.
\/} %\end{lemma}

\medskip
\noindent
Now we are in a position to prove the main results.

\medskip
\noindent
{\it Proof of Theorem 1.1.\/} %\begin{proof}[of Theorem 1.1]

Let $\alpha = 1$, $\mu = -1$, $k = \modul{\kappa} = -\kappa$, and $R > 0$ as in Proposition 3.3.
It is sufficient to prove the assertion for the auxiliary Hamiltonian
$$
  \tilde H_a = -i\sigma_2 {d \over dr} + \sigma_3 + \left({\kappa \over r}
  + {a \over r^2}\chi_{(0, R)}(r)\right) \sigma_1 + {c \over r}
$$
instead of $H_a$, since the eigenvalues of $H_a$ are within $\modul{a}/R$ of
those of $\tilde H_a$.

Introducing the Pr\"ufer transformation
$$
  u = \modul{u} \pmatrix{\sin\thet \cr -\cos\thet \cr}
$$
in the eigenvalue equation $(\tilde H_a - \lambda)\,u = 0$, we find the Pr\"ufer equation
for $\thet$
$$
  \thet' = -\cos 2\thet - \left({\kappa \over r} + {a \over r^2}\chi_{(0, R)}(r)\right)\,\sin 2\thet
         + {c \over r} - \lambda,
$$
which in view of the above choices for $\alpha, \mu$ and $k$ coincides with (\ref{eq:31}) on $(0, R)$ and
with (\ref{eq:32}) on $(R, \infty)$.

For the Pr\"ufer radius $\modul{u}$ we have
$$
  (\log \modul{u})'(r) = - \sin 2\thet + \left({\kappa \over r} + {a \over r^2}\chi_{(0, R)}(r)\right)\,
                \cos 2\thet.
$$
Hence, if $\thet = \Theta_0$ on $(0, R)$, where $\Theta_0$ is the distinguished solution from Lemma
3.1 a) (with $\alpha = 1$), we find
$(\log \modul{u})'(r) \sim \modul{a}/r^2$ $(r\rightarrow 0)$, and hence
$$
  \int_0 \modul{u}^2(r)\,dr \sim {\rm const} \int_0 e^{-2\modul{a}/r}\,dr
  < \infty.
$$
Thus $\Theta_0$ is the Pr\"ufer angle of an $L_2(0, \Any)$ solution of the eigenvalue equation for
$\tilde H_a$.

Similarly, for $a = 0$ we find for a solution $u$ with Pr\"ufer angle $\thet = X_0$, where
$X_0$ is the distinguished solution from Lemma 3.1 b),
$$
  (\log \modul{u})'(r) \sim {\sqrt{k^2 - c^2} \over r},
$$
so
$\modul{u}(r) \sim \mathop{\rm const} r^{\sqrt{k^2 - c^2}}$
$(r \rightarrow 0)$,
whereas for all other solutions $v$
$$
  (\log \modul{v})'(r) \sim {-\sqrt{k^2 - c^2} \over r}
$$
and hence
$\modul{v}(r) \sim \mathop{\rm const} r^{-\sqrt{k^2 - c^2}}$
$(r \rightarrow 0)$.
Thus $X_0$ is the Pr\"ufer angle of the principal solution of the eigenvalue equation for $H_0$.

Analogously,  if (for either $a < 0$ or $a = 0$) $\thet = X_\infty$ on $(R, \infty)$, where
$X_\infty$ is the distinguished solution from Lemma \ref{lem:34}, we have
$$
  (\log \modul{u})'(r) = -\sin 2 X_\infty - {k \over r} \cos 2 X_\infty
  \sim -\sqrt{1 - \lambda^2}
\qquad
  (r \rightarrow \infty),
$$
so
$$
  \int^\infty \modul{u}^2(r)\,dr
  \sim {\rm const} \int^\infty e^{-2\sqrt{1 - \lambda^2}r}\,dr < \infty.
$$
Thus $X_\infty$ is the Pr\"ufer angle of an $L^2(\Any, \infty)$ solution of the eigenvalue
equations for $H_0$ and $\tilde H_a$.

As a consequence, the eigenvalues of $\tilde H_a$ are the (isolated) values of $\lambda$ at
which the
monotone decreasing continuous function $\Theta_0(R, a, c, \Any)$ and the
monotone increasing continuous function $X_\infty(R, c, \Any)$ take the same
value mod $\pi$.
Similarly, the eigenvalues of $H_0$ are the intersection points $\mod \pi$ of
$X_0(R, c, \Any)$ and $X_\infty(R, c, \Any)$.
Hence the assertion follows in view of the uniform convergence of
$\Theta_0(R, a, c, \Any)$ to $X_0(R, c, \Any)$ as $a \rightarrow 0$ (Prop.
3.3).
\hfill $\square$ %\qed\end{proof}

\medskip
\noindent
{\it Proof of Theorem 1.2.\/} %\begin{proof}[of Theorem 1.2]
Let $\alpha = -1$, $\mu = 1$, $k = \modul{\kappa} = \kappa$, and $R > 0$ as in Proposition
3.3. As in the preceding proof, it is sufficient to show the assertion for the auxiliary
Hamiltonian $\tilde H_a$.
We now use the Pr\"ufer transformation
$$
  u = \modul{u} \pmatrix{\cos\thet \cr \sin\thet \cr},
$$
which leads to the Pr\"ufer equation
$$
  \thet' = \cos 2\thet + \left({\kappa \over r} + {a \over r^2}\chi_{(0, R)}(r)\right)\,\sin 2\thet
         + {c \over r} - \lambda,
$$
which in view of the above choices for $\alpha, \mu$ and $k$ coincides with (\ref{eq:31}) on $(0, R)$ and
with (\ref{eq:32}) on $(R, \infty)$.
By studying the asymptotics of $\modul{u}$ as above, we again find that $\Theta_0$ from
Lemma 3.1 a) corresponds to
an $L^2(0, \Any)$ solution of the eigenvalue equation for $\tilde H_a$, $X_0$ from Lemma 3.1 b)
corresponds to the principal solution of the eigenvalue equation for $H_0$, and $X_\infty$
corresponds to an $L^2(\Any, \infty)$ solution of either eigenvalue equation.

Hence the assertion follows as in the preceding proof.
\hfill $\square$ %\qed\end{proof}

\section{Appendix}
\label{sec:app}

In the proofs of this paper we frequently use the following fundamental
observations about first-order ordinary differential equations.
The first lemma (cf. \cite{H} p. 27) is sometimes called \v Caplygin's inequality,
but actually goes back to Peano (\cite{H} p. 44).

\medskip
\noindent
{\bf Lemma 4.1} %\begin{lemma}
\label{lem:a1}
{\it
Let $I \subset \R$ be an interval, $x_0 \in I$ and
$f_j : I \times \R \rightarrow \R$ locally integrable in the first, and
locally Lipschitz continuous in the second argument, $j \in \{1, 2\}$,
with
$f_1(x, y) \le f_2(x, y)$ $(x \in I, y \in \R)$.
Furthermore, let $y_1^0 \le y_2^0$, and $y_j$ be the solution of the
initial value problem
$$
  y'(x) = f_j(x, y),
\qquad
  y(x_0) = y_j^0
\qquad
  (j \in \{1, 2\}).
$$
Then $y_1(x) \le y_2(x)$ $(x \in I, x \ge x_0)$.
\/} %\end{lemma}

\medskip
\noindent
An immediate consequence is the following stability criterion.

\medskip
\noindent
{\bf Lemma 4.2} %\begin{lemma}
\label{lem:a2}
{\it
Let $I \subset \R$ be an interval, $f : I \times \R \rightarrow \R$
locally integrable in the first, and locally Lipschitz continuous in the second argument.
The interval $[y_1, y_2]$ is stable for the differential equation
\begin{equation}
\label{eq:a1}
  y'(x) = f(x, y)
\end{equation}
on $I$ if $f(x, y_1) > 0$, $f(x, y_2) < 0$ $(x \in I)$.
\/} %\end{lemma}

\medskip
\noindent
Here an interval $J$ is called {\it stable\/} on $I$ for (\ref{eq:a1}) if
$y(x_0) \in J \Rightarrow y(x) \in J$ $(x \in I, x \ge x_0)$
for all $x_0 \in I$.

We now use these observations to prove Lemma 2.1.

\medskip
\noindent
{\it Proof of Lemma 2.1.\/} %\begin{proof}[of Lemma 2.1]

\noindent
a)
Let $c \in (-k, 0)$.
Define $J_\pm := \{\thet \in \R \mathrel{|} \pm(c + k \sin 2\thet) > 0\}$, and
$j_\pm^{(\eps)} := \{\thet \in J_\pm \mathrel{|} \modul{\thet - \thet_\pm(c) + j\pi} > \eps
(j \in \Z)\}$ $(\eps > 0)$.

Then for each $\varepsilon > 0$ there
are $\gamma > 0$, $P_0 > \ro_+(c)$ such that for $\varrho \ge P_0$
$$
  \pm(c + (k + \alpha/\varrho) \sin \vartheta) > \gamma
\qquad (\vartheta \in J_\pm^{(\eps)}).
$$

Now consider a solution $\vartheta(\Any, c)$ of (\ref{eq:21}).
If for all $\varepsilon > 0$ there is $P > 0$ and $j \in \Z$ such that
$$
  \modul{\vartheta(\varrho, c) - \varrho_-(c)} \le \varepsilon
\qquad
  (\varrho > P),
$$
this means that
$\lim\limits_{\varrho\rightarrow\infty} \vartheta(\varrho, c)
 = \varrho_-(c) - j \pi$ for some $j \in \Z$.

Otherwise, let $\varepsilon$ be already so small that this is not true,
and $\gamma, P_0 > 0$ as above.
Then there is $\varrho_0 > P_0$ such that
$\modul{\vartheta(\varrho_0, c) - \vartheta_-(c) + j\pi} > \varepsilon$
$(j \in \Z)$, and hence
$$
  \thet(\ro_0, c) \in J_+^{(\eps)} \cup J_-^{(\eps)}
   \cup \bigcup_{j \in \Z} [\thet_+(c) - j\pi - \eps, \thet_+(c) - j\pi + \eps].
$$
In $J_\pm^{(\eps)}$ we have $\pm\vartheta'(\varrho, c) > \gamma/\varrho$,
so if $\vartheta(\varrho_0, c) \in J_\pm^{(\eps)}$, then
$$
  \pm\vartheta(\varrho, c)
  \ge \pm\vartheta(\varrho_0, c) - \gamma \log {\varrho / \varrho_0}
$$
as long as $\thet(\ro, c)$ remains in $J_\pm^{(\eps)}$;
consequently there is $\ro_1$ and $j \in \Z$ such that
$\thet(\ro_1, c) \in [\thet_+(c) - j\pi - \eps, \thet_+(c) - j\pi + \eps]$.

The latter interval is stable by Lemma 4.2.
As $\eps > 0$ was arbitrary, it follows that
$\lim\limits_{\ro\rightarrow\infty} \thet(\ro, c) = \thet_+(c) \mod \pi$.

\smallskip
\noindent
b)
Let $\Gamma$ be an open interval with $\overline{\Gamma} \subset (-k, 0)$,
and $\hat\ro > \ro_+(c)$ $(c \in \Gamma)$.
For $\hat\thet \in \R$ and $c \in \Gamma$,
denote by $\thet(\ro, c, \hat\thet)$ the solution of (\ref{eq:21}) with
initial value $\thet(\hat\ro, c, \hat\thet) = \hat\thet$,
and consider the sets
$$
  S_j := \{(c, \hat\thet) \in (-k, 0) \times \R \mathrel{|}
  \lim_{\ro\rightarrow\infty} \thet(\ro, c, \hat\thet)
   = \thet_+(c) - j\pi\},
$$
$j \in \Z$.
$S_j$ is open. Indeed, let $(c_0, \thet_0) \in S_j$, $\delta := (\thet_+(c_0 - \pi/4)/2$.
Since $\thet_+(c)$ depends continuously on $c$, there is $\eps > 0$
such that
$$
  [\thet_+(c) - \delta, \thet_+(c) + \delta] \subset (\pi/4, \pi)
\qquad
  (c \in [c_0 - \eps, c_0 + \eps]).
$$
Now let $\ro_0 \ge \hat\ro$ be so large that
$\modul{\thet(\ro_0, c_0, \thet_0) - \thet_+(c_0) + j\pi} < \delta/2$.
As $\thet(\ro_0, c, \hat\thet)$ depends continuously on
$(c, \hat\thet)$,
there is positive $\tilde\eps < \eps$ with
$$
  \modul{\thet(\ro_0, c, \hat\thet) - \thet_+(c_0) + j\pi} < \delta
\qquad
  (c \in [c_0 - \tilde\eps, c_0 + \tilde\eps],
  \hat\thet \in [\thet_0 - \tilde\eps, \thet_0 + \tilde\eps]),
$$
and therefore
$\thet(\ro_0, c, \hat\thet) + j\pi$ lies in the stable interval $(\pi/4, \pi)$.

By a) $\lim\limits_{\ro\rightarrow\infty} \thet(\ro, c, \hat\thet)
 = \thet_+(c) - j\pi$,
i.e.
$$
  (c, \hat\thet) \in S_j
\qquad
  (c \in [c_0 - \tilde\eps, c_0 + \tilde\eps],
  \hat\thet \in [\thet_0 - \tilde\eps, \thet_0 + \tilde\eps]).
$$

\smallskip
\noindent
c)
For each $c \in (-k, 0)$, there is exactly one solution
$\thet_\infty(\Any, c)$ such that
$\lim\limits_{\ro \rightarrow\infty} \thet_\infty(\ro, c) = \thet_-(c)$.
Indeed, there is at least one, as $c$ is in some suitable set $\Gamma$
and the corresponding sets
$$
  \Sigma_j(c) = \{\hat\thet \in (0, \pi/4) \mathrel{|}
   (c, \hat\thet) \in S_j\},
$$
$j \in \{0, 1\}$, are nonempty by a) and open by b).

Assume there are two solutions $\thet_1 < \thet_2$ with
$\lim\limits_{\ro\rightarrow\infty} \thet_i(\ro) = \thet_-(c)$,
$i \in \{1, 2\}$.
Then $\ph := 2 (\thet_2 - \thet_1) > 0$ satisfies
$$
\begin{array}{rl}
  \ro\,\ph'(\ro)
  &= 2 (k + \alpha/\ro) (\sin 2 \thet_2 - \sin 2 \thet_1)
\cr\cr
  &= 2 (k + \alpha/\ro) (\sin 2 \thet_1\,{\cos\ph - 1 \over \ph}
    + \cos 2 \thet_1\,{\sin\ph \over \ph})\,\ph
\cr\cr
  &\sim 2 \sqrt{k^2 - c^2}\, \ph > 0
\qquad
  (\ro \rightarrow \infty),
\cr
\end{array}
$$
contradicting $\ph(\ro) \rightarrow 0$.

\smallskip
\noindent
d)
Finally, for fixed $\hat\ro > 0$, $\thet_\infty(\hat\ro, c)$ is strictly
monotone decreasing:
if $c_1 < c_2$ and $\thet_\infty(\hat\ro, c_1) \le \thet_\infty(\hat\ro, c_2)$,
then by Lemma 4.1
$$
  \thet_-(c_1)
  = \lim_{\ro\rightarrow\infty} \thet_\infty(\ro, c_1)
  \le \lim_{\ro\rightarrow\infty} \thet_\infty(\ro, c_2)
  = \thet_-(c_2),
$$
which is not true.

Also, $\thet_\infty(\hat\ro, \Any)$ is continuous.
Indeed, if $c_n\nearrow\hat c$ $(n\to\infty)$, then
$\thet_\infty(\hat\ro, c_n) > \thet_\infty(\hat\ro, \hat c)$
because of the monotonicity, so
$\lim\limits_{n\to\infty} \thet_\infty(\hat\ro, c_n)
 \ge \thet_\infty(\hat\ro, \hat c)$.
However, `$>$' would imply
$(\hat c, \lim \limits_{n\to\infty} \thet_\infty(\hat\ro, c_n)) \in S_0$
in contradiction to the facts that $S_0$ is open and
$(c_n, \thet_\infty(\hat\ro, c_n)) \notin S_0$ $(n \in \N)$.
The right continuity follows in the same way.
\hfill $\square$ %\qed\end{proof}

\medskip
\noindent
{\it Acknowledgement.\/}
H.K. wishes to thank E.B.~Davies for his hospitality at King's College London during the summer of 2002 when the present work was begun.


\begin{thebibliography}{}

\bibitem{AKS}
Arnold V, Kalf H, Schneider A: 
Separated Dirac operators and asymptotically constant linear systems. 
Math. Proc. Cambridge Philos. Soc. \textbf{121} (1997) no. 1, 141--146

\bibitem{B1}
Behncke H:
The Dirac equation with an anomalous magnetic moment.
Math. Z. \textbf{174} (1980), 213--225

\bibitem{B2}
Behncke H:
The Dirac equation with an anomalous magnetic moment II.
in: \textit{Proceedings of the 7th Conference on Ordinary and Partial
Differential Equations, Lect. Notes in Math.} \textbf{964}
(Springer, Berlin 1982)

\bibitem{B3}
Behncke H:
Spectral properties of the Dirac equation with anomalous magnetic moment.
J. Math. Phys. \textbf{26} (1985) no. 10, 2556--2559

\bibitem{GST}
Gesztesy F, Simon B, Thaller B:
On the self-adjointness of Dirac operators with anomalous magnetic moment.
Proc. AMS \textbf{94} (1985), 115--118

\bibitem{H}
Hartman Ph:
\textit{Ordinary differential equations.}
(J. Wiley \& Sons, New York 1964)

\bibitem{K}
Kalf H:
A limit-point criterion for separated Dirac operators and a little known result on Riccati's equation. 
Math. Z. \textbf{129} (1972), 75--82

\bibitem{P1}
Pauli W: Relativistic Field Theories of Elementary Particles. Reviews Mod. Phys. \textbf{13},
(1941) 203--232

\bibitem{P2}
Pauli W:
\textit{Wiss. Briefwechsel mit Bohr, Einstein, Heisenberg u.a. Bd.II, 1930-1939.}
Karl von Meyenn ed. (Springer, Berlin 1985)

\bibitem{P3}
Pauli W:
\textit{Wiss. Briefwechsel mit Bohr, Einstein, Heisenberg u.a. Bd.III, 1940-1949.}
Karl von Meyenn ed. (Springer, Berlin 1993)

\bibitem{RS}
Reed M, Simon B:
\textit{Methods of modern mathematical physics I: Functional analysis.}
(Academic Press, New York 1980)

\bibitem{Th}
Thaller B:
\textit{The Dirac equation.}
(Springer, Berlin 1992)

\end{thebibliography}
\end{document}